\documentclass[12pt]{article}
\usepackage{amsmath,amsfonts,amssymb}
\usepackage[margin=2cm]{geometry}
\usepackage{enumitem}
\usepackage{hyperref}
\usepackage{authblk}
\usepackage[section]{placeins}
\usepackage{todonotes}
\setuptodonotes{inline}

\newtheorem{theorem}{Theorem}[section]

\DeclareMathOperator{\PENT}{PENT}
\DeclareMathOperator{\PLS}{PLS}
\DeclareMathOperator{\STS}{STS}
\DeclareMathOperator{\PBD}{PBD}

\title{Automorphisms of pentagonal geometries\\ with block size 3}
\author{Grahame Erskine\thanks{\texttt{grahame.erskine@open.ac.uk} (corresponding author)} }
\author{Terry S. Griggs\thanks{\texttt{terry.griggs@open.ac.uk}}}

\affil{School of Mathematics and Statistics\\
		The Open University\\
		Walton Hall\\
		Milton Keynes MK7 6AA\\
		UNITED KINGDOM\\}
\date{}

\begin{document}
\maketitle
\begin{abstract}
\noindent A \emph{pentagonal geometry} $\PENT(k,r)$ is a $k$-uniform, $r$-regular partial linear space with the property that the set of points not collinear with a given point $x$ is itself a line in the geometry. The existence spectrum of these geometries is known; we begin the study of automorphisms of these structures in the case $k=3$. We determine the spectrum of pentagonal geometries $\PENT(3,r)$ with either a cyclic or reverse automorphism to be $r \equiv 1$ or $10 \pmod {12}$. In addition, we show that there are 584 pairwise non-isomorphic $\PENT(3,7)$s and classify them by the orders of their automorphism groups.
 \end{abstract}
\noindent\textbf{MSC2020 classification:} 05B25.\\
\noindent\textbf{Keywords:} pentagonal geometry, cyclic automorphism, reverse automorphism.

\section{Introduction}\label{sec:intro}
Pentagonal geometries were introduced in \cite{BBDS}. A \emph{partial linear space} is an ordered pair $(V,\mathcal{L})$ where $V$ is a set of elements, usually called \emph{points}, of cardinality $v$ and $\mathcal{L}$ is a family of subsets of $V$, usually called \emph{lines} or \emph{blocks}, such that every pair of distinct points is contained in at most one line. The number of lines is denoted by $b$. If every line has the same cardinality $k \geq 2$, the space is said to be \emph{uniform}, and if every point is incident with the same number $r \geq 1$ of lines, it is said to be \emph{regular}. Denote such a uniform, regular partial linear space by $\PLS(k,r)$. A \emph{pentagonal geometry}, $\PENT(k,r)$, is a partial linear space $\PLS(k,r)$ in which for all points $x \in V$, the points not collinear with $x$ are themselves collinear. We call this line the \emph{opposite line} to $x$ and denote it by $x^{opp}$. If two points $x$ and $y$ have the same opposite line $x^{opp}=y^{opp}=l$, then $z^{opp}=m$ for all points $z \in l$ where $m$ is the line joining $x$ and $y$. Similarly $w^{opp}=l$ for all $w \in m$. We call such a pair of lines $\{l,m\}$ an \emph{opposite line pair}.

The basic theory of pentagonal geometries $\PENT(k,r)$ is developed in \cite{BBDS}. In this paper we will be concerned with the case $k=3$. From the theory, in a $\PENT(3,r)$, the number of points $v=2r+4$ and the number of lines $b=2r(r+2)/3$. Thus a necessary condition for existence is that $r \equiv 0$ or $1 \pmod 3$. The following result was proved in \cite{GS}.

\begin{theorem}\label{thm:exist}
The existence spectrum of pentagonal geometries $\PENT(3,r)$ is $r \equiv 0$ or $1 \pmod 3$ except $r=4$ or $6$.
\end{theorem}

The pentagonal geometry $\PENT(3,1)$ consists of two disjoint lines forming an opposite line pair, called a \emph{degenerate} pentagonal geometry and $\PENT(3,3)$ is the Desargues configuration. Both are the unique pentagonal geometries for these parameters. Further results on pentagonal geometries with block size 3 are contained in \cite{F} and \cite{FGS}.

The prime objective of this paper is to begin the study of automorphisms of pentagonal geometries $\PENT(3,r)$. Such a geometry is said to be \emph{cyclic} if it admits an automorphism of order $v$ and \emph{reverse} if it admits an automorphism consisting of $v/2$ 2-cycles. These are standard terms applied to all combinatorial structures with the proviso that if $v$ is odd, a reverse automorphism consists of $(v-1)/2$ 2-cycles and a single fixed point. We give constructions for and determine the existence spectrum of both cyclic and reverse pentagonal geometries with block size 3. Finally, we also present a complete enumeration of pentagonal geometries $\PENT(3,7)$ together with details of their automorphism groups.
\section{Cyclic systems}\label{sec:cyclic}
Let $(V,\mathcal{L})$ be a cyclic pentagonal geometry $\PENT(3,r)$. We can assume that $V=\mathbb{Z}_v$ and the automorphism is generated by the mapping $i \mapsto i+1 \pmod v$. We first develop a simple necessary condition for the existence of such a geometry.

Suppose that $r=3s,~s \geq 1$. Then $v=6s+4$ and $b=s(6s+4)$. Under the automorphism, all orbits of 3-sets are full length and thus we require $s$ of these in order to construct the geometry. The opposite line to the point $0,~0^{opp}=\{a,v/2,v-a\}$ for some $a:~1 \leq a \leq v-1$. This is contained in the orbit generated by the block $\{0,v/2-a,v-2a\}$ which contains repeated pairs of points and so cannot be used in the construction of the geometry. So there is no cyclic $\PENT(3,r)$ for $r \equiv 0 \pmod 3$.

Now suppose that $r=3s+1,~s \geq 0$. Then $v=6s+6$ and $b=(3s+1)(6s+6)/3$. In this case, in addition to the orbits of full length, there is a single short orbit of length $v/3$ generated by the block $\{0,v/3,2v/3\}$. Applying the same analysis as in the above case, if $a=v/6$ this orbit will contain all the opposite lines and in particular $0^{opp}=\{v/6,v/2, 5v/6\}=\{s+1,3(s+1),5(s+1)\}$. Since $v$ is even, we can think of a residue mod $v$ as odd or even. Note that the points of $0^{opp}$ are either all odd or all even depending on whether $s$ is even or odd respectively. Let $D$ be the set of odd residues; then the cardinality of $D$ is $v/2=3s+3$. Pairs $\{0,x\},~x \in D$ are contained in the geometry if and only if $x \notin 0^{opp}$. We count these; let the number be $T$.

\begin{tabular}{l}
If $s=4t$, then $r=12t+1,~v=24t+6$ and $T=12t+3-3=12t$.\\
If $s=4t+1$, then $r=12t+4,~v=24t+12$ and $T=12t+6$.\\
If $s=4t+2$, then $r=12t+7,~v=24t+18$ and $T=12t+9-3=12t+6$.\\
If $s=4t+3$, then $r=12t+10,~v=24t+24$ and $T=12t+12$.\\
\end{tabular}

No pairs $\{0,x\},x \in D$ are contained in the short orbit and in a full orbit there are precisely no such pairs or exactly four pairs. Thus there can be no cyclic $\PENT(3,r)$ for $r \equiv 4$ or $7 \pmod {12}$ and we are just left with $r \equiv 1$ or $10 \pmod {12}$ as a necessary condition.

The above analysis not only gives a necessary condition for existence of a cyclic $\PENT(3,r)$, it also establishes its structure: the opposite lines are precisely the blocks of the short orbit. The condition is also sufficient, for which we can appeal to a paper by Wang \& Chang \cite{WC}.

A $(gu,g,k,\lambda)$-\emph{difference family} in $\mathbb{Z}_{gu}$ is a family $F$ of $k$-subsets of $\mathbb{Z}_{gu}$ with the property that the list of differences $\Delta F = \bigcup_{B \in F}\Delta B$ is $\lambda$ times $\mathbb{Z}_{gu} \setminus \{0,u,2u,\ldots,(g-1)u\}$ where $\Delta B = \{b_i-b_j: 1 \leq i,j \leq k, i \neq j\}$ if $B = \{b_1,b_2,\ldots,b_k\}$. A \emph{group divisible design} with block size $k$ and of index $\lambda$, denoted by $(k,\lambda)$-GDD, is a triple $(V,\mathcal{G},\mathcal{B})$ where $V$ is a set of cardinality $v$, $\mathcal{G}$ is a partition of $V$ into parts called \emph{groups} and $\mathcal{B}$ is a family of $k$-subsets of $V$ called \emph{blocks} such that every pair of distinct elements of $V$ occurs in either precisely $\lambda$ blocks or one group but not both. If the partition of $V$ is into groups of cardinality $g$, then the $(k,\lambda)$-GDD is said to be of \emph{type} $g^u$. The existence of a $(gu,g,k,\lambda)$-difference family is equivalent to that of a cyclic $(k,\lambda)$-GDD of type $g^u$ without short orbits. In what follows we will only be concerned with the case where $\lambda = 1$ and will use the more usual notation of $3$-GDD.  

A cyclic $\PENT(3,r)$ is then obtained from a cyclic $3$-GDD of type $6^u$ where $u \equiv 0$ or $1 \pmod 4$ by adjoining the short orbit generated by the block $\{0,2u,4u\}$. That such a difference family exists is stated as part of \cite[Theorem 3.11]{WC}, the proof relying on an earlier result of Jiang \cite{J}.

Alternatively, and probably more directly, cyclic pentagonal geometries with block size $3$ can be constructed from a particular type of integer sequence. A \emph{(pure) Skolem sequence} of order $n$ is a partition of the set $\{1,2,\ldots,2n\}$ into ordered pairs $(a_r,b_r),~r= 1,2,\ldots,n$ such that $b_r-a_r =r$. Such a sequence exists if and only if $n \equiv 0$ or $1 \pmod 4$ \cite{S}. A \emph{hooked Skolem sequence} of order $n$ is a similar partition of the set $\{1,2,\ldots,2n+1\} \setminus \{2n\}$ and exists if and only if $n \equiv 2$ or $3 \pmod 4$ \cite{OK}. Given either a pure or hooked Skolem sequence, it is well known that the collection of blocks  $\{0,n+a_r,n+b_r\}$ generate a cyclic Steiner triple system of order $6n+1$. A \emph{split Skolem sequence} of order $n$ is a similar partition of the set $\{1,2,\ldots,2n+1\} \setminus \{n+1\}$ and exists if and only if $n \equiv 0$ or $3 \pmod 4$ \cite{R1}. Constructions of split Skolem sequences are also given in \cite[p.34]{CR}. In a similar but not quite identical way as for cyclic Steiner triple systems, for $r \equiv 1$ or $10 \pmod {12}$ a split Skolem sequence of order $n=(r-1)/3$ can be used to generate a cyclic $\PENT(3,r)$. Full orbits are generated by the blocks $\{0,n+1+a_r,n+1+b_r\}$ and the short orbit by the block $\{0,v/3,2v/3\}$ under the action of the mapping $i \mapsto i+1 \pmod v$.

The construction is illustrated in Table~\ref{tab:skolem} for the first four cases.

\begin{table}[ht]
    \centering
    \begin{tabular}{|ccc|p{5.5cm}|p{8cm}|}
        \hline
        $r$ & $v$ & $n$ & Split Skolem sequence & Orbit generators\\
        \hline
        10 & 24 & 3 & 
        $(1,3),(2,5),(6,7)$ & 
        $\{0,5,7\},\{0,6,9\},\{0,10,11\},\{0,8,16\}$ \\[2ex]
        13 & 30 & 4 &
        $(1,2),(3,6),(4,8),(7,9)$ &
        $\{0,6,7\},\{0,8,11\},\{0,9,13\},\{0,12,14\},$
        
        $\{0,10,20\}$\\[2ex]
        22 & 48 & 7 &
        $(1,7),(2,6),(3,5),(4,11),$
        
        $(10,15),(9,12),(13,14)$ &

        $\{0,9,15\},\{0,10,14\},\{0,11,13\},\{0,12,19\},$
        
        $\{0,18,23\},\{0,17,20\},\{0,21,22\},\{0,16,32\}$ \\[2ex]
        25 & 54 & 8 &
        $(1,8),(2,5),(10,16),(12,14),$
        
        $(3,4),(6,11),(7,15),(13,17)$ &
        $\{0,10,17\},\{0,11,14\},\{0,19,25\},\{0,21,23\},$
        
        $\{0,12,13\},\{0,15,20\},\{0,16,24\},\{0,22,26\},$
        
        $\{0,18,36\}$ \\
       \hline
    \end{tabular}
    \caption{Construction of cyclic $\PENT(3,r)$ via split Skolem sequences}
    \label{tab:skolem}
\end{table}





\vspace*{1ex}

We have the following result.
\begin{theorem}\label{thm:cycexist}
The existence spectrum of cyclic pentagonal geometries $\PENT(3,r)$ is $r \equiv 1$ or $10 \pmod {12}$.
\end{theorem}

Finally in this section we have some computer generated enumeration results. The number of pairwise non-isomorphic cyclic $\PENT(3,r)$ for the first five values of $r$ are given in Table~\ref{tab:enumer}.

\begin{table}[ht]\centering
\begin{tabular}{|ccccc|}
\hline
$r$ &  & $v$ &  & Number\\
\hline
1 & & 6 & & 1\\
10 & & 24 & & 9\\
13 & & 30 & & 20\\
22 & & 48 & & 20,464\\
25 & & 54 & & 298,880\\
\hline
\end{tabular}
\caption{Numbers of non-isomorphic cyclic $\PENT(3,r)$}
\label{tab:enumer}
\end{table}

Starter blocks for cyclic $\PENT(3,10)$ and $\PENT(3,13)$ on the set $\mathbb{Z}_v$ under the mapping $i \mapsto i+1 \pmod v$ are given in Table~\ref{tab:starters}. Commas and set brackets are omitted for brevity and the integers 10 to 30 are represented by the letters $a$ to $t$ respectively.

\begin{table}[ht]
    \centering
    \begin{tabular}{|p{4cm}|p{5cm}|}
    \hline
    $\PENT(3,10)$  &  $\PENT(3,13)$ \\
    \hline
    \texttt{013 05e 06h 08g} & \texttt{013 04d 06m 07i 0ak} \\
    \texttt{013 05e 06d 08g} & \texttt{013 04c 06j 07g 0ak} \\
    \texttt{013 05f 06d 08g} & \texttt{013 04d 06m 07j 0ak} \\
    \texttt{016 02h 03d 08g} & \texttt{013 04c 06h 07l 0ak} \\
    \texttt{016 029 03d 08g} & \texttt{013 04l 06e 07j 0ak} \\
    \texttt{016 029 03e 08g} & \texttt{013 04c 06j 07l 0ak} \\
    \texttt{01b 027 039 08g} & \texttt{013 04l 06e 07i 0ak} \\
    \texttt{01b 027 03i 08g} & \texttt{013 04c 06h 07g 0ak} \\
    \texttt{01b 02j 039 08g} & \texttt{017 02b 03h 04c 0ak} \\
     & \texttt{017 02b 03h 04m 0ak} \\
     & \texttt{017 02b 03g 04c 0ak} \\
     & \texttt{017 02b 03g 04m 0ak} \\
     & \texttt{017 02l 03h 04c 0ak} \\
     & \texttt{017 02l 03h 04m 0ak} \\
     & \texttt{017 02l 03g 04c 0ak} \\
     & \texttt{017 02l 03g 04m 0ak} \\
     & \texttt{018 02j 039 04g 0ak} \\
     & \texttt{018 02j 039 04i 0ak} \\
     & \texttt{018 02d 039 04g 0ak} \\
     & \texttt{018 02d 03o 04g 0ak} \\
    \hline
    \end{tabular}
    \caption{Starter blocks for cyclic $\PENT(3,10)$ and $\PENT(3,13)$ systems}
    \label{tab:starters}
\end{table}
\section{Reverse systems}\label{sec:reverse}
Let $(V,\mathcal{L})$ be a reverse pentagonal geometry $\PENT(3,r)$. Denote the automorphism by $I$. First note that a cyclic pentagonal geometry is also a reverse pentagonal geometry with $I$ as the mapping $i \mapsto i+v/2 \pmod v$. Hence there exist reverse pentagonal geometries $\PENT(3,r)$ for $r \equiv 1$ or $10 \pmod {12}$. As in the case of cyclic systems, we show that this condition is also necessary.

Let $x \in V$. Note that $I(x)\in x^{opp}$ since if $\{x,I(x),z\}$ were a line in the geometry then $z$ would be a fixed point of $I$ which is impossible. Thus $\ell=x^{opp}=\{I(x),y,z\}$ for some $y,z \in V$. This implies that $m=y^{opp}=\{I(y),x,w\},~w \in V$ and by applying the automorphism that $I(y)^{opp}=\{y,I(x),I(w)\}$ and hence that $w=I(z)$. It follows that $z^{opp}=y^{opp}=m$ and therefore the lines $\ell$ and $m$ are an opposite line pair. So the opposite lines are partitioned under $I$ into opposite line pairs and $v=2r+4$ is divisible by $6$, eliminating the possibility that $r \equiv 0 \pmod 3$.

Continuing with the analysis, adjoin three new points $\infty_0$, $\infty_1$, $\infty_2$ and a new line $\{\infty_0, \infty_1, \infty_2\}$ to the pentagonal geometry. For each opposite line pair $\{I(x),y,z\}$, $\{x,I(y),I(z)\}$ also adjoin lines $\{\infty_0,x,I(x)\}$, $\{\infty_0,y,I(y)\}$, $\{\infty_0,z,I(z)\}$, $\{\infty_1,y,I(z)\}$, $\{\infty_2,z,I(y)\}$, $\{\infty_1,z,x)\}$, $\{\infty_2,I(x),I(z)\}$, $\{\infty_1,I(x),I(y)\}$, $\{\infty_2,y,x\}$.
We have a Steiner triple system on $2r+7$ points, $\STS(2r+7)$, which has a reverse automorphism under the involution $I$ and $(\infty_1, \infty_2)$ with $\infty_0$ as a fixed point. These exist if and only if $2r+7 \equiv 1,3,9$ or $19 \pmod{24}$, \cite {D,R2,T}. Hence $r \equiv 1$ or $10 \pmod{12}$ since the case $r \equiv 0 \pmod 3$ has already been eliminated.

\newpage
As for cyclic systems we have the following result.
\begin{theorem}\label{thm:revexist}
The existence spectrum of reverse pentagonal geometries $\PENT(3,r)$ is $r \equiv 1$ or $10 \pmod {12}$.
\end{theorem}

Finally, we present two further constructions of reverse $\PENT(3,r)$.
The first of these begins with a reverse Steiner triple system, $\STS(v)$, on base set $\{\infty,x,I(x):0 \leq x \leq (v-3)/2\}$ where the reverse automorphism $I$ is the set of cycles $(x,I(x))$ and $\infty$ is the fixed point. Delete all the blocks containing the point $\infty$ to obtain a $3$-GDD of type $2^{(v-1)/2}$. This also has a reverse automorphism. Weight each point by a factor $3$ and replace each point $x$ or $I(x)$ by a block $\{x_0,x_1,x_2\}$ or $\{I(x)_0,I(x)_1,I(x)_2\}$ respectively. These will be the opposite line pairs of the pentagonal geometry. For each pair of blocks $\{x,y,z\}$, $\{I(x),I(y),I(z)\}$ of the GDD replace the former by the blocks
$\{x_0,y_0,z_0\}$, $\{x_1,y_1,z_1\}$, $\{x_2,y_2,z_2\}$, $\{x_0,y_1,z_2\}$,
$\{x_1,y_2,z_0\}$, $\{x_2,y_0,z_1\}$, $\{x_0,y_2,z_1\}$, $\{x_1,y_0,z_2\}$,
$\{x_2,y_1,z_0\}$ and the latter by applying the automorphism $I$. This construction gives reverse $\PENT(3,r)$ for $r \equiv 1,10,25$ or $34 \pmod{36}$.

The second construction is similar and uses the concept of a \emph{pairwise balanced design}, $\PBD(v,K,\lambda)$. This is a pair $(V,\mathcal{B})$ where $V$ is a set of cardinality $v$ and $\mathcal{B}$ is a family of subsets of $V$ called \emph{blocks} with cardinalities in $K$ and such that every pair of distinct elements of $V$ occurs in precisely $\lambda$ blocks. We will also require two reverse $3$-GDDs, one of type $6^4$ and the other of type $6^5$. These are given below.

The construction begins with a $\PBD(v,\{4,5\},1)$. These exist for $v \equiv 0$ or $1\pmod 4$, $v \notin \{8,9,12\}$, see~\cite[p.250]{HB}. Weight each point of the set $V$ by a factor 6, i.e. replace each point $x$ by six points $(x,j),~0 \leq j \leq 5$. Create lines $\{(x,0),(x,1),(x,2)\}$ and $\{(x,3),(x,4),(x,5)\}, x \in V$. These will be the opposite line pairs of the pentagonal geometry. Now replace each block, say $\{w,x,y,z\}$, of cardinality 4 in the pairwise balanced design by the reverse $3$-GDD of type $6^4$ assigning the four points $w,x,y,z$ to the first co-ordinate $0,1,2,3$ in any order and identifying the second co-ordinates. 

Proceed similarly for the blocks of cardinality 5. By construction, each of the two reverse GDDs below has an involutory automorphism of the form $(x,i)\mapsto(x,i+3)$ and this automorphism is preserved in the $\PENT(3,r)$; this construction thus gives reverse $\PENT(3,r)$ for $r \equiv 1$ or $10 \pmod{12}, r \notin \{22,25,34\}$.  

\begin{tabular}{l}
\textbf{Reverse~$3$-GDD of type $6^4$}\\
Base set $V= \{(i,j):0 \leq i \leq 3,~0 \leq j \leq 5\}$.\\
Partition $G$ is into groups $G_i =\{(i,j):0 \leq j \leq5\},~0 \leq i \leq 3$.\\
Blocks. There are 72 blocks, 9 of which are the following\\
$\{(0,0),(1,0),(2,0)\}, \{(0,1),(1,1),(2,1)\}, \{(0,2),(1,2),(2,2)\},$\\
$\{(0,0),(1,1),(2,2)\}, \{(0,1),(1,2),(2,0)\}, \{(0,2),(1,0),(2,1)\},$\\
$\{(0,0),(1,2),(2,1)\}, \{(0,1),(1,0),(2,2)\}, \{(0,2),(1,1),(2,0)\}.$\\
A further 27 blocks are obtained by applying the mappings\\
(i) $(0,j) \mapsto (0,j+3),~(1,j) \mapsto (3,j)$;\\
(ii) $(1,j) \mapsto (1,j+3),~(2,j) \mapsto (3,j)$;\\
(iii) $(2,j) \mapsto (2,j+3),~(0,j) \mapsto (3,j)$ to the above 9 blocks.\\
Finally the remaining 36 blocks are obtained by applying the mapping $(i,j) \mapsto (i,j+3)$\\
to all of the above blocks, arithmetic in the second co-ordinate modulo 6. 
\end{tabular}

\begin{tabular}{l}
\textbf{Reverse~$3$-GDD of type $6^5$}\\
Base set $V= \{(i,j):0 \leq i \leq 4,~0 \leq j \leq 5\}$.\\
Partition $G$ is into groups $G_i =\{(i,j):0 \leq j \leq5\},~0 \leq i \leq 4$.\\
Blocks. There are 120 blocks, 12 of which are the following\\
$\{(0,0),(1,0),(2,3)\}, \{(0,1),(1,1),(2,5)\}, \{(0,2),(1,2),(2,4)\},$\\
$\{(0,0),(1,1),(2,4)\}, \{(0,1),(1,2),(2,3)\}, \{(0,2),(1,0),(2,5)\},$\\
$\{(0,0),(1,2),(3,2)\}, \{(0,1),(1,0),(3,0)\}, \{(0,2),(1,1),(3,1)\}.$\\
$\{(0,0),(1,4),(3,1)\}, \{(0,1),(1,3),(3,2)\}, \{(0,2),(1,5),(3,0)\}.$\\
A further 48 blocks are obtained by applying the mappings\\
$(i,j) \mapsto (i+k,j),~k \in\{1,2,3,4\}$ to the above 12 blocks, \\
with arithmetic in the first co-ordinate modulo 5.\\
Finally the remaining 60 blocks are obtained by applying the mapping $(i,j) \mapsto (i,j+3)$\\
to all of the above blocks, arithmetic in the second co-ordinate modulo 6.
\end{tabular}

In general, the reverse pentagonal geometries obtained from the two constructions will not be cyclic.  We thus have reverse, non-cyclic $\PENT(3,r)$ for all admissible values of $r$ except for $r=22$. Such a pentagonal geometry is given below. It has 48 points and 352 lines and we present it as a set of orbit starters under a cyclic subgroup of automorphisms of order 24 generated by the permutation:

\begin{tabular}{rl}
    $a=$ & \texttt{(0 n a f 5 j 6 e 1 l b g 3 k 7 c 2 m 9 h 4 i 8 d)} \\
    & \texttt{(o K z D s I u B q J y F r H w A p L x E t G v C)}
\end{tabular}

Orbit starters:

\texttt{012 06j 07E 08G 0cA 0eL 0fs 0gw 0pH 0qB 0rz 0tv 0uF 0xJ opq ouH}

Note that $a^{12}$ is the involution \texttt{(0 3)(1 4)(2 5)}\ldots and so the orbits of length 8 generated by \texttt{012} and \texttt{opq} give the opposite line pairs; all other orbits have length 24.
\section{PENT(3,7) systems}\label{sec:pent37} 
We have already observed that the pentagonal geometries $\PENT(3,r)$ are unique for $r=1$ or $3$ and that there are no such geometries for $r=4$ or $6$. The next case to consider is $r=7$. In~\cite[Lemma 2.4]{BBDS}, it was proved that pentagonal geometries $\PENT(k,2k+1)$ either have no opposite line pairs or the points are partitioned into opposite line pairs. Further, in~\cite[Lemma 10.2]{BBDS}, with the aid of a computer check it was shown that for $\PENT(3,7)$, the first possibility does not exist. For the latter possibility the $18$ points are partitioned into three opposite line pairs $\{l_0,m_0\},\{l_1,m_1\},\{l_2,m_2\}$ with the remaining 36 lines coming from a group divisible design, $3$-GDD of type $6^3$, (i.e. a Latin square of side $6$), whose groups are the points of $l_i \cup m_i,~i= 0,1,2$.

There are precisely $12$ paratopy or main classes of Latin square of side 6 which are listed in~\cite[p.137]{HB}. Thus the pentagonal geometries $\PENT(3,7)$ can be classified into $12$ different types corresponding to the main classes of Latin square, \emph{not} isomorphism classes as stated incorrectly in \cite{GS}.  There are 584 pairwise non-isomorphic pentagonal geometries $\PENT(3,7)$ partitioned by main class of Latin square in the order in which they appear in \cite{HB} and order of automorphism group. They are given in Table~\ref{tab:pent37}.

\begin{table}[ht]\centering
\begin{tabular}{|c|cccccccccccc|c|}
\hline
Aut & \multicolumn{12}{c|}{Main class} & \\
group & 1 & 2 & 3 & 4 & 5 & 6 & 7 & 8 & 9 & 10 & 11 & 12 & TOTAL\\
\hline
1 &  & 1 & 1 & 45 & 95 & 5 &  & 8 & 26 & 32 & 30 &  & 243\\
2 & 1 & 7 & 4 & 29 & 50 & 12 &  & 18 & 24 & 15 & 18 & 1 & 179\\
3 &  & 1 &  &  &  &  &  & 1 &  & 3 & 3 &  & 8\\
4 & 4 & 5 & 4 & 11 & 18 & 10 & 1 & 10 & 13 & 2 & 4 & 2 & 84\\
6 & 3 & 1 &  &  &  &  & 1 & 2 &  & 4 & 2 & 2 & 15\\
8 & 1 & 4 & 1 & 2 & 4 & 2 & 2 & 4 & 3 &  & 2 & 1 & 26\\
12 &  &  &  &  &  & 1 &  &  &  &  &  &  & 1\\
16 & 1 & 2 &  &  &  &  & 1 & 2 &  &  &  &  & 6\\
24 & 2 & 1 & 1 &  &  & 1 &  &  & 1 &  & 2 & 2 & 10\\
48 & 2 & 1 &  &  &  &  & 2 & 2 &  &  &  &  & 7\\
72 &  &  &  &  &  & 1 &  &  &  &  &  &  & 1\\
144 &  & 1 &  &  &  &  &  &  &  &  &  &  & 1\\
432 & 1 &  &  &  &  &  &  &  &  &  &  &  & 1\\
648 &  &  &  &  &  &  &  &  &  &  &  & 1 & 1\\
1296 &  &  &  &  &  &  & 1 &  &  &  &  &  & 1\\
\hline
TOTAL & 15 & 24 & 11 & 87 & 167 & 32 & 8 & 47 & 67 & 56 & 61 & 9 & 584\\
\hline
\end{tabular}
\caption{Numbers of $\PENT(3,7)$ by main class of Latin square}
\label{tab:pent37}
\end{table}

From the results of Sections~\ref{sec:cyclic} and~\ref{sec:reverse} above, none of these systems can have a cyclic or reverse automorphism so all automorphisms of order $2$ will have fixed points. However, there are four systems which are transitive. These are those with the largest automorphism groups but none admits a transitive subgroup of order $18$.

\begin{flushleft}
	\bibliographystyle{abbrv}
	\bibliography{PENT}
\end{flushleft}

\section*{Statements and declarations}
The authors have no financial or other conflicts of interest to declare.

Data files produced during this research are available from the corresponding author on request.

\end{document}